\documentclass[11pt,twoside,a4paper]{article}
\usepackage{amsfonts}
\usepackage{amssymb}
\usepackage{amsmath}

\begin{document}

\newcommand{\pa}{\partial}
\newcommand{\opa}{\overline\pa}
\newcommand{\ol}{\overline}

\numberwithin{equation}{section}

\newcommand\C{\mathbb{C}}  
\newcommand\R{\mathbb{R}}
\newcommand\Z{\mathbb{Z}}
\newcommand\N{\mathbb{N}}
\newcommand\PP{\mathbb{P}}

{\LARGE \centerline{On holomorphic foliations admitting invariant CR manifolds}}
\vspace{0.8cm}

\centerline{\textsc {\large Judith Brinkschulte}\footnote{Universit\"at Leipzig, Mathematisches Institut, PF 100920, D-04009 Leipzig, Germany. 
E-mail: brinkschulte@math.uni-leipzig.de
 }}

\vspace{0.5cm}

\begin{abstract} 
We study holomorphic foliations of codimension $k\geq 1$ on a complex manifold $X$ of dimension $n+k$ from the point of view of the exceptional minimal set conjecture. For $n\geq 2$ we show in particular that if the holomorphic normal bundle $N_{\mathcal{F}}$ is Griffiths positive, then the foliation does not admit a compact  invariant set that is a  complete intersection of $k$ smooth real hypersurfaces in $X$.
\end{abstract}

\vspace{0.5cm}

\section{Introduction}

Let $X$ be a complex manifold and $\mathcal{F}$ a (singular) holomorphic foliation of codimension $k\geq 1$ on $X$. In complex dynamics, one is interested in understanding the structure of the set of accumulation point of the leaves of $\mathcal{F}$. A natural question asks under what condition on $X$ and/or $\mathcal{F}$
does it hold that every leaf of $\mathcal{F}$ accumulates to the singular set $\mathrm{Sing}(\mathcal{F})$?\\

In this context, Brunella in \cite{Br} stated the following conjecture, which may also be formulated for holomorphic foliations of codimension $k\geq 1$:\\
Let $X$ be a compact connected complex manifold of dimension $\geq 3$, and let $\mathcal{F}$ be a codimension one holomorphic foliation on $X$ whose normal bundle $N_{\mathcal{F}} = TX/T\mathcal{F}$ is ample. Then 
every leaf of $\mathcal{F}$ accumulates to  $\mathrm{Sing}(\mathcal{F})$.\\

Assume to the contrary that not every leaf of $\mathcal{F}$ accumulates to  $\mathrm{Sing}(\mathcal{F})$: Then $X$ contains a nonempty compact subset $\mathcal{M}$ which is invariant by $\mathcal{F}$ and disjoint from $\mathrm{Sing}(\mathcal{F})$ (a so-called exceptional minimal set). The aim of this paper is to prove that such an $\mathcal{M}$ cannot be a smooth intersection of $k$ real hypersurfaces whose normals are complex linearly independent.\\

More precisely, our main result is the following\\

\pagebreak

\newtheorem{main}{Theorem}[section]
\begin{main}   \label{main}   \ \\
Let $X$ be a complex manifold of dimension $n+k$, and let $M\subset X$ be a compact smooth real submanifold that is a complete intersection of $k$ real hypersurfaces in $X$.  Suppose that on some neighborhood $U$ of $M$, there exists a codimension $k$ holomorphic foliation $\mathcal{F}$ which leaves $M$ invariant. Then, for $n\geq 2$, the holomorphic normal bundle $N_{\mathcal{F}}$ of $\mathcal{F}$ in $U$ does not admit any Hermitian metric $h$ such that $(N_{\mathcal{F}},h)$ is Griffiths positive  on $U$.
\end{main}

By a complete intersection we mean that $M$ is defined by $k$ smooth function $\rho_j\in\mathcal{C}^\infty(X)$:
$$M = \lbrace z\in X\mid \rho_1(z)= \ldots =\rho_k(z)=0\rbrace,$$
with $\pa\rho_1\wedge\ldots\wedge \pa\rho_k \not=0$ on $M$.\\

Theorem \ref{main} generalizes the main result of \cite{Br} in two ways. First of all, we pass from codimension one to codimension $k\geq 1$. Second, $X$ is an arbitrary complex manifold in Theorem \ref{main}, possibly noncompact, whereas in \cite{Br} it was assumed to be compact K\"ahler.\\

Another interpretation of Theorem \ref{main} is from the point of view of classifying compact Levi-flat CR manifolds. Let us remind that
a CR manifold of type $(n,k)$ is given by a smooth real manifold $M$ of dimension $2n+k$ and a complex subbundle $T^{1,0}M$ of rank $n$ of $\C\otimes M$ that is stable under the Lie-bracket. It is also required that $T^{1,0}M\cap T^{0,1}M = \lbrace 0\rbrace$, where
 $T^{0,1}M = \overline{T^{1,0}M}$.
If moreover $T^{1,0} M + T^{0,1}M$ is closed  under the Lie bracket, then $M$ is called {\it Levi-flat}. It follows from the theorems of Frobenius and Newlander-Nirenberg that Levi-flat CR manifolds are locally foliated by complex $n$-dimensional submanifolds.\\

The normal bundle of a Levi-flat CR manifold of type $(n,k)$  is the  vector bundle of rank $k$ over $M$ defined by
$$N_M = \C\otimes TM / (T^{1,0}M + T^{0,1}M).$$
$N_M$ carries the structure of a CR vector bundle over $M$, that is the restriction to each leaf of the foliation of $M$ is a holomorphic vector bundle.\\

Remark that due to a theorem of Andreotti-Fredricks \cite{AF}, a real-analytic CR manifold always admits a generic CR embedding into a complex manifold. If $M$ is also Levi-flat, then the Levi-foliation of $M$ can be holomorphically extended to some neighborhood $U$ of $M$ (see \cite{R}). Also a metric on $N_M$ with leafwise positive curvature can always be extended to a metric on the extended bundle with positive curvature. Thus Theorem \ref{main} enables us to state the following corollary:

\newtheorem{realana}[main]{Corollary}
\begin{realana}   \label{realana}   \ \\
Let $X$ be a complex manifold of dimension $n+k$, and let $M\subset X$ be a real-analytic compact  CR manifold of type $(n,k)$. If $M$ is Levi-flat and a complete intersection of $k$ real hypersurfaces in $X$, then, for $n\geq 2$, the holomorphic normal bundle $N_M$ of $M$ in $X$ does not admit any Hermitian metric $h$ such that $(N_M ,h)$ is Griffiths positive along the leaves of $M$.
\end{realana}

As a very special case, Corollary \ref{realana} therefore includes the nonexistence results for Levi-flat real hypersurfaces in $\C\PP^n,\ n\geq 3$, because the normal bundle of a real hypersurface in $\C\PP^n$ naturally inherits positivity from the ambient complex manifold (see \cite{LN} and \cite{S}). For $k> 1$, the corresponding nonexistence result for smooth Levi-flat CR manifolds in complex projective spaces can be found in \cite{B1}. We also refer the interested reader to \cite{FSW}, where the authors weaken the assumptions on the Levi-flat CR manifold under a more restrictive assumption on $k$.\\

In \cite{B} Corollary \ref{realana} was proved for $k=1$, also if $M$ is only smooth instead of real-analytic ($M$ is then a compact Levi-flat real hypersurface), thereby generalizing previous results by Brunella \cite{Br} and Ohsawa \cite{O7}.\\

Note that Theorem \ref{main} does not hold for Levi-flat CR manifolds of CR-dimension $n=1$. Counterexamples can be found in \cite{Br} and \cite{O8}.\\

For other results related to the exceptional minimal set conjecture for holomorphic foliations of arbitrary codimension, we refer the reader to \cite{CF}.\\

The organization of this paper is as follows:
We argue by contradiction and assume that we are given a holomorphic foliation of codimension $k$ with $N_{\mathcal{F}}$ Griffiths positive that leaves invariant a compact complete intersection $M$ of $k$ 
real hypersurfaces as in the statement of Theorem \ref{main}. In section 3 we prove that $M$ admits a neighborhood basis having certain $q$-convexity/concavity properties. This allows us to holomorphically extend CR sections of $(\mathrm{det}N_{\mathcal{F}})^{\otimes\ell}$ in section 4. 
Combining a version of Kodaira's embedding theorem for Levi-flat CR manifolds proved by Ohsawa with the holomorphic extension of CR sections of $(\mathrm{det}N_{\mathcal{F}})^{\otimes\ell}$, we then show that a tubular neighborhood of $M$ can be generically embedded into a projective K\"ahler space $\hat X$ (section 5). 
In section 6, we prove some $L^2$-Hodge-type symmetry results on the regular part of $\hat X$. 
The final argument is then given in section 7: we first extend the Chern curvature form of $\mathrm{det}N_{\mathcal{F}}$ to a $d$-exakt $(1,1)$-form on the regular part of $\hat X$. Using the $L^2$-Hodge-type symmetry results from section 6, we may then find a smooth potential for the Chern curvature form of $\mathrm{det}N_{\mathcal{F}}$ on an open neighborhood of $M$. This will contradict the maximum principle on the leaves of $M$ and therefore prove our main result.\\

{\it Remark:} As a referee pointed out, in Theorem \ref{main} and Corollary \ref{realana}, it is possible to relax the assumption on  $M$ of being a complete intersection of $k$ hypersurfaces to the CR dimension of $M$ being equal to $n$ and the line bundle $\mathrm{det}(N_M)$ being topologically trivial. In this situation, in the proof of Proposition \ref{function}, one has to replace the function $\rho$ with the distance function to $M$ with respect to the metric $\omega_o$.

\section{Preliminaries}

Let $Y$ be a complex manifold of dimension $n$ endowed with a Hermitian metric $\omega$, and let $E$ be a holomorphic vector bundle on $Y$ with a Hermitian metric $h$. 
We recall that $(E,h)$ is said to be {\it Griffiths positive} if for all $y\in Y$ and all non-zero decomposable tensors $v\otimes e\in TY\times Y_y$, the curvature term $\langle i\Theta_h(E)(v,v)e,e\rangle$ is positive, where $i\Theta_h(E)\in \mathcal{C}^\infty_{1,1}(Y,\mathrm{Herm}(E,E))$ is the curvature of the Chern connection of $(E,h)$.\\

For integers $0\leq p,q\leq n$, we use the following notations:\\
$\mathcal{D}^{p,q}(Y,E)$ denotes the space of smooth, compactly supported $E$-valued $(p,q)$-forms on $Y$.\\
$L^2_{p,q}(Y,E,\omega,h)$ denotes the Hilbert space obtained by completing $\mathcal{D}^{p,q}(Y,E)$ with respect to the $L^2$-norm $\Vert \cdot\Vert_{\omega,h}$ induced by $\omega$ and $h$.\\

As usual, the differential operator $\opa$ is extended as a densely defined closed linear operator on $L^2_{p,q}(Y,E,\omega,h)$, whose domain of definition is 
$$\mathrm{Dom}\opa = \lbrace f\in L^2_{p,q}(Y,E,\omega,h)\mid \opa f \in L^2_{p,q+1}(Y,E,\omega,h)\rbrace,$$
where $\opa f$ is computed in the sense of distributions. The Hilbert space adjoint of $\opa$ will be denoted by $\opa^\ast \ (= \opa^\ast_{\omega,h})$.\\

We also define the space of harmonic forms,
$$\mathcal{H}^{p,q}(Y,E,\omega,h)= L^2_{p,q}(Y,E,\omega,h)\cap\mathrm{Ker}\opa\cap\mathrm{Ker}\opa^\ast_{\omega,h},$$
 and the $L^2$-Dolbeault cohomology groups of $Y$,
$$H^{p,q}_{L^2}(Y,E,\omega,h) = L^2_{p,q}(Y,E,\omega,h)\cap\mathrm{Ker}\opa / L^2_{p,q}(Y,E,\omega,h)\cap\mathrm{Im}\opa.$$

Whenever we feel that it is clear from the context, we will omit the dependency of the $L^2$-spaces, norms, operators etc. on the hermitian metric $h$ of the vector bundles under considerations. \\

In section 4, we shall also use the following variant of the $\opa$-operator: by $\opa_c$ we denote the strong minimal realization of $\opa$ on $L^2_{p,q}(Y,E,\omega,h)$. This means that $u\in\mathrm{Dom}\opa_c\subset L^2_{p,q}(Y,E,\omega,h)$ if there exists $f\in L^2_{p,q+1}(Y,E,\omega,h)$ and a sequence $(u_\nu)_{\nu\in\N} \subset\mathcal{D}^{p,q}(Y,E)$ such that $u_\nu\longrightarrow u$ and $\opa u_\nu\longrightarrow f = \opa_c u$ in $L^2_{p,q+1}(Y,E,\omega,h)$.\\

The Hilbert space adjoint of $\opa_c$ will be denoted by $\vartheta$; it is the weak maximal realization of the formal adjoint of $\opa$ on $L^2_{p,q}(Y,E,\omega,h)$.

\section{Convexity properties of tubular neighborhoods}

From now on, $M$ will always denote a smooth submanifold of real codimension $k$ in a complex manifold $X$ of dimension $n+k$ that is a complete intersection of $k$ smooth real hypersurfaces $\Sigma_j = \lbrace z\in U\mid \rho_j(z)=0\rbrace$ and that
is invariant by a holomorphic foliation $\mathcal{F}$ on some neighborhood of $M$.\\

The aim of this section is to find a real-valued function on $U\setminus M$, with a certain growth order near $M$, whose Levi-form has $n+1$ positive and $k-1$ negative eigenvalues.

For that purpose, we fix a Hermitian metric $\omega_o$ on $X$ and define $\rho = (\rho_1^2 + \ldots + \rho_k^2)^\frac{1}{2}$.\\

Note the analogy with the following model case: If $X$ is a complex manifold of dimension $n+k$ and $Y\subset X$ is a compact complex submanifold of dimension $n$ in $X$ such that the normal bundle $N_{Y/X}$ of $Y$ in $X$ is positive in the sense of Griffiths, then $Y$ admits a tubular neighborhood basis  consisting of domains $U$ such that the Levi form of $\pa U$ has $n$ negative and $k-1$ positive eigenvalues (see \cite{G}).\\

For $k=1$, the following Proposition was proved in \cite{Br}.\\

\newtheorem{function}{Proposition}[section]
\begin{function}   \label{function}   \ \\
Let $X$ be a complex manifold of dimension $n+k$, and let $M\subset X$ be a smooth compact real submanifold that is a complete intersection of $k$ real hypersurfaces in $X$.  Suppose that on some neighborhood $U$ of $M$, there exists a codimension $k$ holomorphic foliation $\mathcal{F}$ which leaves $M$ invariant. We moreover assume that  the holomorphic normal bundle $N_{\mathcal{F}}$ of $\mathcal{F}$ is positive in the sense of Griffiths on $U$. Then, after possibly shrinking $U$, there exist a smooth nonnegative real-valued function $v$ on $U$ and constants $c, c^\prime >0$ such that if $\lambda_1\leq \ldots\leq \lambda_{n+k}$ denote the eigenvalues of $i\pa\opa(-\log v)$ with respect to $\omega_o$, then 
\begin{enumerate}
\item[(i)] $\rho^2\lambda_j \in \lbrack -c^\prime, -c\rbrack,\quad  j=1,\ldots, k-1,$
\item[(ii)] $\lambda_j\in \lbrack c,c^\prime\rbrack,\quad j = k,\ldots, n+k-1,$
\item[(iii)] $\rho^2\lambda_{n+k} \in \lbrack c,c^\prime\rbrack$
\end{enumerate}
Moreover we have $v= O(\rho^2)$.
\end{function}

{\it Proof.} The proof is an adaptation of the argument given in \cite{O7}.\\

We may choose a finite covering of $M$ by open sets $U_\alpha$, such that $\mathcal{F}$ is defined by $k$ holomorphic 1-forms $\varpi_\alpha^1,\ldots,\varpi_\alpha^k$  on $U_\alpha$ (i.e. $\mathrm{Ker}\varpi_\alpha^j\subset T^{1,0}\mathcal{F}$). Then there exist holomorphic $k\times k$-matrices $G_{\alpha\beta}$ on $U_\alpha\cap U_\beta$ satisfying
\begin{equation} \label{matrices1}
\varpi_\alpha = G_{\alpha\beta}\varpi_\beta,
\end{equation}
with $\varpi_\alpha = (\varpi_\alpha^1,\ldots,\varpi_\alpha^k)^T$; in fact the matrices $G_{\alpha\beta}$ are the transition matrices for the holomorphic normal bundle $N_\mathcal{F}$ with respect to the dual frames $(\varpi_\alpha)^\ast$ of $\varpi_\alpha$.\\

The fiber metric of $N_\mathcal{F}$ may then be identified with a system of positive definite Hermitian matrices $H_\alpha$ on $U_\alpha$ such that 

$$H_\alpha = \big( \langle (\varpi_\alpha^j)^\ast, (\varpi_\alpha^\ell)^\ast\rangle \big)_{j,\ell = 1}^k.$$

We then have
\begin{equation}  \label{matrices2}
H_\alpha = G_{\beta\alpha}^T H_\beta \ol G_{\beta\alpha}^\ast \quad\mathrm{on}\ U_\alpha\cap U_\beta.
\end{equation}

 Since $M$ is a smooth Levi-flat CR manifold in $X$, we may assume that  $M$ is locally defined by $M\cap U_\alpha = \lbrace z\in U_\alpha\mid (\mathrm{Im}f_\alpha^j)(z) = 0, \ j = 1,\ldots, k\rbrace$, where $f_\alpha^j$ is a function such that $\opa f_\alpha^j$ vanishes to infinite order on $M\cap U_\alpha$ and $df_\alpha^1\wedge\ldots\wedge df_\alpha^k\not=0$ on $M\cap U_\alpha$ (this can be seen using the parametric equations for $M$, see e.g. \cite{HT}). \\

On $U_\alpha\cap M$ we have $\varpi_\alpha =  A_\alpha df_\alpha$ for some smooth matrix $A_\alpha$ which is invertible on $U_\alpha$ and holomorphic along the leaves of $M$. From (\ref{matrices1}) and (\ref{matrices2}) it follows that we have
\begin{equation*}
 (A_\alpha df_\alpha)^T  \ol H_\alpha \ol A_\alpha d\ol f_\alpha   = (A_\beta df_\beta)^T \ol H_\alpha \ol A_\beta d\ol f_\beta \ \mathrm{on}\ U_\alpha\cap U_\beta \cap M
\end{equation*}

But then 
\begin{equation*}
(A_\alpha\mathrm{Im}f_\alpha)^T \ol H_\alpha \ol A_\alpha \mathrm{Im}f_\alpha - (A_\beta\mathrm{Im}f_\beta)^T \ol H_\beta \ol A_\beta \mathrm{Im}f_\beta = \mathcal{O}(\rho^3) \ \mathrm{along}\ U_\alpha\cap U_\beta \cap M
\end{equation*}

Therefore, invoking Whitney's extension theorem, there exists a nonnegative real-valued $\mathcal{C}^2$ function $v$ defined in a tubular neighborhood of $M$, smooth away from $M$, such that $v = (A_\alpha\mathrm{Im}f_\alpha)^T \ol H_\alpha \ol A_\alpha \mathrm{Im}f_\alpha + O(\rho^3)$ on $U_\alpha$. \\

It now remains to estimate $-i\pa\opa \log v = -i \frac{\pa\opa v}{v} + i\pa\log v\wedge \opa\log v$ on $U\setminus M$.\\

 It suffices to estimate the Levi-form of $-\log\big( (A_\alpha\mathrm{Im}f_\alpha)^T \ol H_\alpha \ol A_\alpha \mathrm{Im}f_\alpha\big) $ on $U_\alpha \setminus M$. Setting $H_\alpha^\prime = \ol A_\alpha^T H_\alpha A_\alpha$, we have to estimate the Levi-form of $-\log \big(\mathrm{Im} f_\alpha^T \ol H_\alpha^\prime \mathrm{Im}f_\alpha\big)$ on $U_\alpha\setminus M$. \\

Therefore we let $V\in T^{1,0}U$ be a unitary vector that we decompose orthogonally into $V = V_t + V_n$, with $V_t\in \bigcap_{j=1}^k\mathrm{Ker}(\pa\mathrm{Im}f_\alpha^j)$. A straightforward computation then gives that up to terms vanishing  on $M$
\begin{align*}
-i\pa\opa\log v (V_t, \ol V_t)  = & i \Big( \frac{\mathrm{Im}f_\alpha^T\opa\pa \ol H_\alpha^\prime \mathrm{Im}f_\alpha}{v}
 +\frac{\mathrm{Im}f_\alpha^T\pa \ol H_\alpha^\prime \mathrm{Im}f_\alpha \wedge 
\mathrm{Im}f_\alpha^T\opa \ol H_\alpha^\prime \mathrm{Im}f_\alpha}{v^2} \Big) (V_t,\ol V_t)\\
 = & i \Big( \frac{\mathrm{Im}f_\alpha^T \ol H_\alpha^\prime\opa \big( (\ol H_\alpha^\prime)^{-1}\pa \ol H_\alpha^\prime \big)\mathrm{Im}f_\alpha}
{v} \Big) (V_t,\ol V_t) +\\
&  i\Big( \frac{\mathrm{Im}f_\alpha^T\opa \ol H_\alpha^\prime \wedge (\ol H_\alpha^\prime)^{-1}\pa \ol H_\alpha^\prime \mathrm{Im}f_\alpha}{v} \Big) (V_t,\ol V_t)
  + i\frac{\vert \mathrm{Im}f_\alpha^T\pa \ol H_\alpha^\prime \mathrm{Im}f_\alpha \vert^2 (V_t)
}{v^2} .
\end{align*}

Both terms in the last line of the above equation are nonnegative since
\begin{align*}
\Big( \mathrm{Im}f_\alpha^T\opa \ol H_\alpha^\prime \wedge (\ol H_\alpha^\prime)^{-1}\pa \ol H_\alpha^\prime \mathrm{Im}f_\alpha \Big) (V_t, \ol V_t)   = & \ \vert (H_\alpha^\prime)^{-1/2}\opa H_\alpha^\prime \mathrm{Im} f_\alpha (\ol V_t)\vert^2 
\end{align*}

Thus we get 
\begin{equation*}
-i\pa\opa\log v (V_t, \ol V_t) \geq i \Big( \frac{\mathrm{Im}f_\alpha^T \ol H_\alpha^\prime\opa \big( (\ol H_\alpha^\prime)^{-1}\pa \ol H_\alpha^\prime \big)\mathrm{Im}f_\alpha}
{v} \Big) (V_t,\ol V_t).
\end{equation*}

Since $H_\alpha$ represents the metric of $N_{\mathcal{F}}$ in a local holomorphic trivialization over $U_\alpha$ we have

\begin{equation*}
i\Theta(N_{\mathcal{F}}) = i \opa \big( (\ol H_\alpha)^{-1}\pa \ol H_\alpha \big)\ \mathrm{on}\ U_\alpha
\end{equation*}

But as the matrix $A_\alpha$ is holomorphic along the leaves of $M$, we may assume that $\opa A_\alpha$ vanishes to finite order along $M$. But then $i \opa \big( (\ol H_\alpha^\prime)^{-1}\pa \ol H_\alpha^\prime \big)$ is the Chern curvature tensor of $N_{\mathcal{F}}$ in the almost holomorphic trivialization given by the change of basis $A_\alpha$. We may therefore 
find  $c > 0$ such that
\begin{equation} \label{eq1}
-i\pa\opa \log v (V_t,\ol V_t) \geq (c + \epsilon)\omega_o( V_t,\ol V_t),
\end{equation}
where $\epsilon$ can be made as small as we wish by shrinking $U$.\\

Moreover, since $\opa f_\alpha^j$ vanishes to infinite order on $U_\alpha\cap M$ and $V_t\in  \bigcap_{j=1}^k\mathrm{Ker}(\pa\mathrm{Im}f_\alpha^j)$, a careful computation shows that
\begin{equation*}
	\vert i\pa\opa \log v(V_t, \ol V_n)\vert  \lesssim 
	\vert V_t\vert \cdot\vert V_n\vert
\end{equation*}

Therefore we get
\begin{equation}  \label{eq2}
\vert i\pa\opa \log v (V_t,\ol V_n)\vert  \leq \epsilon \omega_o( V_t,\ol V_t) + \frac{C}{\epsilon} \omega_o(V_n,\ol V_n)
\end{equation}
for some constant $C$.\\
 
 The remaining arguments do not involve the Levi-flatness of $M$, therefore we may replace $v$ with $\rho$ for the computations involving $-i\pa\opa\log v$. We have $\pa\rho = \frac{1}{\rho}(\rho_1\pa\rho_1 +\ldots + \rho_k\pa\rho_k)$ and
 \begin{align} \label{hessian}
 	-i\pa\opa\rho = & -\frac{i}{\rho}(\sum_{j=1}^{k} \rho_j\pa\opa\rho_j + \sum_{j=1}^{k} \pa\rho_j\wedge\opa\rho_j) \nonumber\\
 	& + \frac{i}{\rho^3} (\rho_1\pa\rho_1 +\ldots + \rho_k\pa\rho_k) \wedge (\rho_1\opa\rho_1 +\ldots + \rho_k\opa\rho_k)
 \end{align}
 
 The first sum in (\ref{hessian}) gives only a small contribution with respect to the other two terms as the point approaches $M$, and it can therefore be neglected. \\

We now choose a  $(1,0)$ vector field $\xi \notin \bigcap_{j=1}^k \mathrm{Ker}(\pa\rho_j)$ on $U\setminus M$ such that $\pa \rho(\xi) = 1$ on $U\setminus M$ as follows:
$$\xi = \frac{1}{\rho}\sum_{j=1}^k\rho_j \xi_j,$$
where $\xi_j$ are $(1,0)$ vector fields on $U$ satisfying $\pa\rho_\ell (\xi_j) = \delta_{j\ell}$.
A direct computation shows that the sum of the second and third term in the right hand side of (\ref{hessian}) applied to $(\xi,\ol\xi)$ is zero, hence the leading term of
$i\pa\opa\log\rho (\xi,\ol\xi)$ is $i\frac{\pa\rho\wedge\opa\rho}{\rho^2}(\xi,\ol\xi) = \frac{1}{\rho^2}$.\\

We may thus conclude that
\begin{equation} \label{eq3}
-i\pa\opa \log v (\xi,\ol \xi) \geq a\frac{1}{v}
\end{equation}
for some constant $a > 0$ (shrinking $U$ if necessary).\\ Combining (\ref{eq1}), (\ref{eq2}) and (\ref{eq3}), the minimum-maximum theorem shows the existence of $c> 0$ in (ii) and (iii). The existence of $c^\prime > 0$ follows from the reversed inequalities in (\ref{eq1}) (as a consequence of the Levi-flatness of $M$) and (\ref{eq3}) (as follows from the definition of the function $v$).\\

The $(k-1)$ negative eigenvalues in (i) come from the second term on the right hand side of (\ref{hessian}).

\hfill$\square$\\

\section{Holomorphic extension of CR sections of $\mathrm{det}N_M$}

The key result of this section is Proposition \ref{extension}, the extension of CR sections of the normal bundle
to holomorphic sections of $\mathrm{det}N_{\mathcal{F}}$ over $U$. This will enable us to holomorphically embed a tubular neighborhood of $M$ into some complex projective space in the next section. The proof of this holomorphic extension property needs several steps; it is a modification of the arguments in section 7 of \cite{B}.\\

From now on, we will exploit the properties of the function $\varphi = -\log v$, where $v$ as in Proposition \ref{function}. Shrinking $U$ if necessary, we may assume that $v$ is actually defined in an open neighborhood of $\ol U$ and satisfies (i), (ii) and (iii) of Proposition \ref{function} on that open neighborhood.\\

We start with the following lemma that can be proved as Lemma 4.2 in \cite{B2}:

\newtheorem{metric}{Lemma}[section]
\begin{metric}   \label{metric} \ \\
There exists a complete hermitian metric $\omega_{M}$ on $\ol U\setminus M$ with the following properties:
\begin{itemize}
\item[(i)] Let $\gamma_{1} \leq \ldots \leq \gamma_{n+k}$ be the
eigenvalues of $i \pa \opa \varphi$ with respect to $\omega_{M}$. There
exists $\sigma >0$ such that $\gamma_{1}+ \ldots +\gamma_{k} >
\sigma$ on $U\setminus M$.
\item[(ii)] There are constants $a, b > 0$ such that $a \ \omega_o \leq
\omega_{M} \leq b \ \rho^{-2} \omega_o$.
\item[(iii)] There is a constant $C > 0$ such that $\vert \pa \omega_{M}
\vert_{\omega_{M}} \leq C$.
 \end{itemize}
\end{metric}

From Lemma \ref{metric} we obtain the following $L^2$-vanishing result:

\newtheorem{vanish2}[metric]{Proposition}
\begin{vanish2}     \label{vanish2}  \ \\
Let $E\longrightarrow U $ be a holomorphic line bundle. Then there exist  $N\in\N$ such that the following holds: Assume $0\leq q\leq n-1$, and let $f\in L^2_{0,q}(U\setminus M, E, \omega_M, -N\varphi)\cap\mathrm{Ker}\opa$. Then there exists a $(0,q-1)$-form $g\in L^2_{0,q-1}(U\setminus M, E, \omega_M, -N\varphi)$ satisfying $\opa g = f$ and $\Vert g\Vert_{\omega_M, -N\varphi}\leq\Vert f\Vert_{\omega_M,-N\varphi}$.
\end{vanish2}

{\it Proof.} We fix $0\leq q\leq n-1$
It follows from the generalized Bochner-Kodaira-Nakano inequality (see $\cite[\mathrm{Chapter\ VII,\ section\ 1-3}]{D1}$) and standard computations that for $u\in\mathcal{D}^{0,q}(U\setminus M, E)$ we have
\begin{align*}
\frac{3}{2} ( \Vert \opa u\Vert^2_{\omega_M, -N\varphi} + \Vert \opa^\ast u\Vert^2_{\omega_M,-N\varphi})\hspace{10cm}\\
\geq N \ll (\gamma_1 + \ldots + \gamma_{n+k-q})u,u\gg_{\omega_M, -N\varphi}\hspace{7cm}\\ - 
\ll (c_{q+1} + \ldots + c_{n+k})u,u\gg_{\omega_M, -N\varphi} \hspace{7cm}\\
-\frac{1}{2} (\Vert \tau u\Vert^2_{\omega_M, -N\varphi} + \Vert \ol\tau u\Vert^2_{\omega_M, -N\varphi} + \Vert \tau^\ast u\Vert^2_{\omega_M, -N\varphi} + \Vert \ol\tau^\ast u\Vert^2_{\omega_M, -N\varphi}),\hspace{3cm}
\end{align*}
where $c_1\leq \ldots\leq c_{n+k}$ are the  curvature eigenvalues of a fixed hermitian metric of $E$ with respect to $\omega_M$, $\Lambda$ is the adjoint of multiplication by $\omega_M$ and $\tau = \lbrack \Lambda,\pa\omega_M\rbrack$.\\

For $N$ sufficiently big, it therefore follows from the properties of $\omega_M$ described in Lemma \ref{metric} that
\begin{equation} \label{apriori1}
\Vert u\Vert^2_{\omega_M, -N\varphi} \leq \Vert\opa u\Vert^2_{\omega_M, -N\varphi} + \Vert \opa^\ast u\Vert^2_{\omega_M, -N\varphi}
\end{equation}
for $u\in\mathcal{D}^{0,q}(U\setminus M, E)$. By the completeness of $\omega_M$, the estimate (\ref{apriori1}) extends to $u\in L^2_{0,q}(U\setminus M,\omega_M, -N\varphi)\cap\mathrm{Dom}\opa\cap\mathrm{Dom}\opa^\ast$ with compact support in $U$.\\ 

Note that  (\ref{apriori1}) holds also for $q=n$; it is only now that we have to restrict to $q\leq n-1$:
Shrinking $U$ if necessary, we may assume that the boundary of $U$ is smooth and that its Levi form has $n$ negative and $k-1$ positive eigenvalues.
From $\cite[\mathrm{Theorem\ VI\ and\ Theorem\ 7.4}]{G}$ we may therefore deduce that if $N$ is sufficiently big, then
\begin{equation} \label{apriori2}
\Vert u\Vert^2_{\omega_M, -N\varphi} \leq \Vert\opa u\Vert^2_{\omega_M, -N\varphi} + \Vert \opa^\ast u\Vert^2_{\omega_M, -N\varphi}
\end{equation}
for $u\in L^2_{0,q}(U\setminus M,\omega_M, -N\varphi)\cap\mathrm{Dom}\opa\cap\mathrm{Dom}\opa^\ast$ with compact support in $\ol U\setminus K$, where $K$ is a compact containing an open neighborhood of $M$ in $U$. \\

Using two cut-off functions $\chi_1,\chi_2$, where $\chi_1$ has compact support in $U$ and equals one in an open neighborhood of $M$ and $\chi_2 = 1-\chi_1$, we may use (\ref{apriori1}) and (\ref{apriori2}) to conclude that if $N$ is sufficiently big, then
\begin{equation*}
\Vert u\Vert^2_{\omega_M, -N\varphi} \leq \Vert\opa u\Vert^2_{\omega_M, -N\varphi} + \Vert \opa^\ast u\Vert^2_{\omega_M, -N\varphi}
\end{equation*}
for $u\in L^2_{0,q}(U\setminus M,\omega_M, -N\varphi)\cap\mathrm{Dom}\opa\cap\mathrm{Dom}\opa^\ast$, $0\leq q\leq n-1$. From this a priori estimate, the assertion of the Proposition follows in a standard way.
\hfill$\square$\\

In the next section, we want to holomorphically extend CR sections over $M$ of some high tensor power of the line bundle $\mathrm{det}N_{\mathcal{F}}$, which is a positive line bundle by assumption.
On the other hand, we can multiply the metric of $\mathrm{det}N_{\mathcal{F}}$ by $e^{N\varphi}$. This adds $-Ni\pa\opa\varphi$ to the curvature, so the curvature of $\mathrm{det}N_{\mathcal{F}}$ can be made partly negative near on $U\setminus M$ by taking $N$ sufficiently large. This modification, however, would require the CR sections that we wish to extend to be sufficiently regular. 

Since this is false in general,
 we have to 
use some approximation arguments, reducing the involved $\opa$-equation to compactly supported forms. As a result we can prove

\newtheorem{extension}[metric]{Proposition}
\begin{extension}     \label{extension}  \ \\
Let $\ell \in\N$ be sufficiently large, and assume that $s$ is a CR-section of $(\mathrm{det}N_{\mathcal{F}})^{\ell}$ over $M$ of class at least $\mathcal{C}^{k+4}$. Then there exists a holomorphic section $\tilde s$ of $(\mathrm{det} N_{\mathcal{F}})^{\ell}$ on $U$ such that $\tilde s_{\mid M} = s$.
\end{extension}

{\it Proof.} First we choose a $\mathcal{C}^{k+4}$-extension $s_o$ of $s$ to $X$ such that $\opa s_o$ vanishes to the order $k+3$ along $M$, i.e. $\vert\opa s_o \vert^2_{\omega_o} = O(\vert\rho\vert^{2k+6})$. \\

Now we consider an exhaustion of $U\setminus M$ by  domains $V_\varepsilon = \lbrace z\in U\mid v (z) > \varepsilon^2 \rbrace$; we recall that $v = O(\rho^2)$. Moreover, we define the annular domains $D_j = V_{\frac{1}{j}}\setminus \ol V_{\frac{2}{j}}$. \\

Then we choose a sequence of smooth cut-off functions $\chi_j$ with compact support in $V_{\frac{1}{j}}$ such that $\chi_j\equiv 1$ on $\ol V_{\frac{2}{j}}$ and $\vert d\chi_j\vert^2_{\omega_M} \leq 1$ (this is possible since $\omega_M$ is complete on $\ol U\setminus M$). Then 
\begin{equation}  \label{f_j}
f_j := \opa(\chi_j\opa s_o)\in L^2_{0,2}(U\setminus M, (\mathrm{det} N_{\mathcal{F}})^{\ell},\omega_M)\cap\mathrm{Ker}\opa
\end{equation}
 is compactly supported in $D_j$. \\
Applying Lemma \ref{csupport} yields $u_j\in L^2_{0,1}(U\setminus M, (\mathrm{det} N_{\mathcal{F}})^{\ell},\omega_M)$ supported in $D_j$ satisfying $\opa u_j = f_j$ and

\begin{align*}
 \Vert u_j\Vert^2_{\omega_M}  \leq &  C^2 j^4\Vert f_j\Vert^2_{\omega_M} \lesssim j^4\Vert\opa \chi_j\wedge\opa s_o\Vert^2_{\omega_M}\\
\leq & j^4\int_{D_j} \vert\opa s_o\vert^2_{\omega_M} dV_{\omega_M}\leq j^4\int_{D_j}\rho^{-2k}
\vert\opa s_o\vert^2_{\omega_o} dV_{\omega_o} \lesssim 1
\end{align*}

Now $g_j = \chi_j\opa s_o - u_j$ is $\opa$-closed and supported in $V_{\frac{1}{j}}$, hence compactly supported in $\ol U\setminus M$.  By Proposition \ref{vanish2}, there exists $N\in\N$ such that we can find solutions $h_j\in L^2_{0,0}(U\setminus M, (\mathrm{det}N_{\mathcal{F}})^{\ell},\omega_M, -N\varphi)$ satisfying $\opa h_j = g_j$. Hence $g_j\in L^2_{0,1}(U\setminus M, (\mathrm{det}N_{\mathcal{F}})^{\ell},\omega_M,  -k\varphi)\cap\mathrm{Im}\opa$. By Lemma \ref{bottom}, we can therefore have $ h_j\in 
L^2_{0,0}(U\setminus M, (\mathrm{det}N_{\mathcal{F}})^{\ell},\omega_M, -\varphi)$ and  $\Vert  h_j\Vert_{\omega_M,-\varphi}\leq C_o \Vert g_j\Vert_{\omega_M ,-\varphi}$.
But
\begin{align*}
\Vert g_j\Vert^2_{\omega_M,-\varphi} \lesssim &\  \Vert \chi_j\opa s_o\Vert^2_{\omega_M,-\varphi} + \Vert u_j\Vert^2_{\omega_M,-\varphi}\\
 \lesssim & \ \int_{V_{\frac{1}{j}}} \vert \opa s_o\vert^2_{\omega_o} \rho^{-2k-2} dV_{\omega_o} +  \Vert u_j\Vert^2_{\omega_M,-\varphi}\\
 \lesssim &\  \int_{V_{\frac{1}{j}}} \rho^{2k+6} \rho^{-2k-2} dV_{\omega_o} + j^4\int_{D_j}\rho^{-2k-2}
 \vert\opa s_o\vert^2_{\omega_o} dV_{\omega_o}  \lesssim 1
\end{align*}

 Therefore the sequence $( h_j)$ is bounded in $L^2_{0,0}(U\setminus M, (\mathrm{det}N_{\mathcal{F}})^{\ell},\omega_M,  -\varphi)$, hence has a subsequence that weakly converges to\\ $ h_o\in L^2_{0,0}(U\setminus M, (\mathrm{det} N_{\mathcal{F}})^{\ell},\omega_M, -\varphi)$. Since $\opa  h_j = \opa s_o$ on $V_{\frac{2}{j}}$, we must therefore have $\opa h_o = \opa s_o$ in $U\setminus M$. \\
Moreover, since $ h_o\in L^2_{0,0}(U\setminus M, (\mathrm{det} N_{\mathcal{F}})^{\ell},\omega_M,  -\varphi)$, we have $$\int_{U\setminus M}\vert  h_o\vert^2 \rho^{-2} dV_{\omega_M} < +\infty.$$
 This clearly implies that the trivial extension of $ h_o$ to $U$ satisfies $\opa h_o = \opa s_o$ as distributions on  $U$ (not only on $U\setminus M$). Hence $h_o$ is of class at least $\mathcal{C}^{k+3}$ by the hypoellipticity of $\opa$, and must therefore vanish on $M$.\\

Thus $\tilde s = s_o -  h_o$ is a holomorphic section of $(\mathrm{det}N_{\mathcal{F}})^{\ell}$ over $U$ extending $s$. \hfill$\square$\\

\newtheorem{csupport}[metric]{Lemma}
\begin{csupport}     \label{csupport}  \ \\
Let $\ell \in\N$ be sufficiently large and $f_j$ be defined by (\ref{f_j}).
For some constant $C> 0$, independent of $j\in\N$,  there exists $u_j \in L^2_{0,1}(U\setminus M, (\mathrm{det}N_{\mathcal{F}})^{\ell},\omega_M)$, supported in $D_j$, such that $\opa u_j = f_j$ and
$$\Vert u_j\Vert_{ \omega_M} \leq C j^2\Vert f_j\Vert_{\omega_M}.$$
\end{csupport}

{\it Proof.} To abbreviate notations, we set $F = (\mathrm{det}N_{\mathcal{F}})^{\ell}$.\\
Note that the boundary of $D_j$ consists of two parts: the part $\pa V_{\frac{1}{j}}$ whose Levi-form has $n$ positive and $k-1$ negative eigenvalues and the part $-\pa V_{\frac{2}{j}}$ whose Levi-form has $n$ negative and $k-1$ positive eigenvalues. Since the case $k=1$ was treated in \cite{B} we may assume $k\geq 2$. But then $D_j$ satisfies condition $Z(n+k-1)$ (see \cite{FK}), hence the $\opa$-Neumann problem satisfies subelliptic estimates in degree $(p,n+k-1)$ for all $0\leq p\leq n+k$.  In particular this implies (see also \cite{FK}) that
$$\opa : L^2_{n+k, n+k-2}(D_j, F^\ast, \omega_o) \longrightarrow  L^2_{n+k, n+k-1}(D_j, F^\ast, \omega_o)$$
has closed range (subelliptic estimates are proved  using local computations, therefore they are valid independent of the curvature of $F$, since $F$ is defined over $\ol D_j$).
\\

Now we use a duality argument from \cite{Ch-S}: Let $\opa_c$ be the strong minimal realization of $\opa$ on $L^2_{0,1}(D_j, F,\omega_o)$. Then by Theorem 3 of \cite{Ch-S} the range of $\opa_c$ is closed in $L^2_{0,2}(D_j, F, \omega_o)$, and $\opa_c$-exact forms $f\in
L^2_{0,2}(D_j, F, \omega_o)$ are characterized by the usual orthogonality condition:
$$\int_{D_j} f \wedge\eta = 0   \quad\forall \eta \in L^2_{n+k,n+k-2}(D_j,F^\ast,\omega_o)\cap\mathrm{Ker}\opa $$
 But, using Stokes' theorem, we get for $\eta \in \mathcal{C}^\infty_{n+k,n+k-2}(\ol D_j,F^\ast)\cap\mathrm{Ker}\opa$
$$ \int_{D_j} \opa(\chi_j\opa s_o) \wedge\eta = \int_{\pa D_j}(\chi_j\opa s_o)\wedge\eta = -\int_{\pa V_{\frac{2}{j}}} \opa s_o\wedge\eta = -  \int_{\pa V_{\frac{2}{j}}} \opa (s_o\wedge\eta)  = 0,$$
and this also holds for $\eta \in L^2_{n+k,n+k-2}(D_j,F^\ast,\omega_o)\cap\mathrm{Ker}\opa$ using a standard approximation argument and the subelliptic estimates in degree $(n+k,n+k-1)$.

 Hence $f_j = \opa(\chi_j\opa s_o)$ belongs to the image of $\opa_c$, i.e. there exists $u_j\in L^2_{0,1}(D_j, F, \omega_o)$ satisfying $\opa_c u_j = f_j$. As usual, we assume that $u_j$ is the minimal $L^2$-solution i.e. $u_j \in L^2_{0,1}(D_j, F,\omega_o)\cap(\mathrm{Ker}\opa_c)^\perp \subset\mathrm{Ker}\vartheta$.
 In particular, $u_j$ is smooth on $\ol D_j$, and the trivial extension of $u_j$ by zero outside $\ol D_j$ (which we still denote by $u_j$), satisfies $\opa u_j = f_j$ as distributions on $U\setminus M$ (by definition of the strong minimal realization $\opa_c$). It remains to estimate $\Vert u_j\Vert_{\omega_M}$.\\

First we note that it follows from the subelliptic estimates in degree $(n+k,n+k-1)$ that $u_j$ is sufficiently smooth on $\ol D_j$:  $\ast u_j$ is of bidegree $(n+k,n+k-1)$ and satisfies the elliptic system $\opa^\ast(\ast u_j)= \ast f_j$, $\opa(\ast u_j) = 0$. Since $f_j$ is of class $\mathcal{C}^3$ and vanishes outside a compact of $D_j$, $u_j$ is at least in the Sobolev space $W^3$ and smooth up to the boundary outside the support of $f_j$.\\ 
 
We will now estimate $u_j$ by using a priori estimates on the  domains $W_j = U\setminus\ol V_{\frac{2}{j}}$. The Levi-form of $\pa W_j$ has $n$ negative and $k-1$ positive eigenvalues. Since  we may assume $k\geq 2$, it follows in particular that the Levi-form of $\pa W_j$ has at least one positive eigenvalue everywhere.

We now modify the metric in $\mathrm{det}N_{\mathcal{F}}^\ast$ by a bounded factor $\exp{(-m\rho^2)}$. This adds to the curvature a term which is $$mi\pa\opa\rho^2= \sum_{j=1}^k \big( 2m\rho_j i\pa\opa\rho_j + 2mi\pa\rho_j\wedge\opa\rho_j\big).$$
Taking $m$ sufficiently large, and shrinking $U$ if necessary, we may therefore assume that the curvature form of  $\mathrm{det}N_{\mathcal{F}}^\ast$ has at least $k\geq 2$ positive eigenvalues. But then
the curvature term in the Bochner-Kodaira-Nakano identity of $\mathrm{det}N_{\mathcal{F}}^\ast$ is positive in degree $(n+k,n+k-1)$.\\

It then follows from the Bochner-Kodaira-Nakano formula with boundary term 
(see $\cite[\mathrm{Theorem\ 7.2}]{G})$ that there exists $\lambda >0$ such that
$$\Vert v\Vert^2_{\omega_o, W_j} \leq \frac{\lambda }{\ell m} (\Vert \opa v\Vert^2_{\omega_o,W_j} + \Vert\opa^\ast v\Vert^2_{\omega_o, W_j})$$
for all $v\in L^2_{n+k,n+k-1}(W_j,  F^\ast,\omega_o) \cap\mathrm{Dom}\opa\cap\mathrm{Dom}\opa^\ast$.
From this we infer by Serre duality as in \cite{Ch-S} that
\begin{equation}  \label{estimatewj}
\Vert v\Vert^2_{\omega_o, W_j} \leq \frac{\lambda}{\ell m} (\Vert \opa_c v\Vert^2_{\omega_o, W_j} + \Vert \vartheta v\Vert^2_{\omega_o,W_j})
\end{equation}
for all $v\in L^2_{0,1}(W_j, F,\omega_o) \cap\mathrm{Dom}\opa_c\cap\mathrm{Dom}\vartheta$.\\

We now choose an extension $\tilde u_j$ of $u_j$ to $\tilde V_j = V_{\frac{1}{2j}}\setminus \ol V_{\frac{2}{j}}$ such that $\tilde u_j \in \mathrm{Dom}\opa_c\cap\mathrm{Dom}\vartheta$ (on $\tilde V_j$!) and such that $$\Vert\vartheta\tilde u_j\Vert^2_{\omega_o, \tilde V_j} + \Vert \opa\tilde u_j\Vert^2_{\omega_o, \tilde V_j} \leq b( \Vert\vartheta u_j\Vert^2_{\omega_o, \tilde D_j} + \Vert \opa u_j\Vert^2_{\omega_o, D_j} + \Vert u_j\Vert^2_{\omega_o, D_j})$$ for some constant $b$ not depending on $u_j$ nor on $j$. This is possible for $j$ sufficiently large by general Sobolev extension methods (locally we flatten the boundary $\pa D_j$ and extend the sufficiently smooth $u_j$ componentwise  across $\pa D_j$ by first order reflection, then we use a partition of unity).\\

Next we choose a smooth cut-off function $\tilde\chi_j$ with compact support in $ V_{\frac{1}{2j}}\setminus V_{\frac{2}{j}}$  such that $\tilde\chi_j\equiv 1$ on $\ol D_j$. This is possible with an estimate $\vert d\tilde\chi_j\vert^2 \lesssim j^2$.\\

Applying (\ref{estimatewj}) with $\chi_j \tilde u_j$ yields
\begin{align*}
\Vert u_j\Vert^2_{\omega_o, D_j} \leq \Vert \tilde\chi_j \tilde u_j\Vert^2_{\omega_o, W_j} \leq & \frac{\lambda}{\ell m} (\Vert \opa_c (\tilde\chi_j\tilde u_j)\Vert^2_{\omega_o, W_j} + \Vert \vartheta (\tilde\chi_j\tilde u_j)\Vert^2_{\omega_o, W_j}) \\
 \lesssim & \frac{\lambda}{\ell m} (\Vert \opa_c \tilde u_j\Vert^2_{\omega_o, W_j} + \Vert \vartheta \tilde u_j\Vert^2_{\omega_o, W_j} + j^2 \Vert \tilde u_j\Vert^2_{\omega_o, W_j}) \\
 \leq & \frac{\lambda}{\ell m} b ( \Vert\vartheta u_j\Vert^2_{\omega_o, D_j} + \Vert \opa u_j\Vert^2_{\omega_o, D_j} + j^2\Vert u_j\Vert^2_{\omega_o, D_j})\\
 = & \frac{\lambda}{\ell m} b ( \Vert\opa_c u_j\Vert^2_{\omega_o, D_j} + j^2\Vert u_j\Vert^2_{\omega_o, D_j})
\end{align*}
For $m = m(j) \sim j^2$, the additional weight factor $\exp (-m\rho^2)$ is bounded on $D_j$. Hence we may choose first $m$ and then $\ell$ sufficiently large so that
$$\Vert u_j\Vert_{\omega_o} \leq \Vert f_j\Vert_{\omega_o}.$$

It remains to compare the norms $\Vert u_j\Vert_{\omega_o}$ and $\Vert u_j\Vert_{\omega_M}$.  From the construction of $\omega_M$ (see Lemma \ref{metric} ), we have $dV_{\omega_M} \sim \rho^{-2k} dV_{\omega_o}$. We also have $\vert f_j\vert^2_{\omega_o}   \lesssim \rho^{-4} \vert f_j\vert^2_{\omega_M}$ since $f_j$ is a $(0,2)$-form. On the other hand, we have $\omega_M \gtrsim \omega_o$, which implies $\vert u_j\vert^2_{\omega_o} \gtrsim \vert u_j\vert^2_{\omega_M}$. Since $u_j$ is supported in $D_j$, we thus have
$$\Vert u_j\Vert^2_{\omega_M} \lesssim j^{2k} \Vert u_j\Vert^2_{\omega_o} \leq j^{2k} \Vert f_j\Vert^2_{\omega_o} \lesssim j^4 \Vert f_j\Vert^2_{\omega_M},$$
which proves the desired estimate.
\hfill$\square$\\

The point of the following lemma is that even though $\ell\in\N$ can be arbitrary big, the weight function $-\varphi$ stays the same (it does not have to be multiplied by a large integer as $\ell$ increases!).\\

\newtheorem{bottom}[metric]{Lemma}
\begin{bottom}     \label{bottom}  \ \\
Let $\ell\in\N$ be arbitrary. Then there exists a constant $C_o$ such that for every
$u \in L^2_{0,0}(U\setminus M, (\mathrm{det}N_{\mathcal{F}})^{\ell},\omega_M,-\varphi)$ with $\opa u\in L^2_{0,1}(U\setminus M,  (\mathrm{det}N_{\mathcal{F}})^{\ell},\omega_M, -\varphi)$ satisfies  $\Vert u\Vert_{\omega_M,-\varphi} \leq C_o \Vert f\Vert_{\omega_M,-\varphi}$.
\end{bottom}

{\it Proof.}
We start by working with the metric $\omega_o$ on $U$. 
Recall that in degree $(0,0)$, the curvature term in the Bochner-Kodaira-Nakano identity is given by minus the trace of the curvature form with respect to the metric under consideration.
Therefore it is convenient to have 
\begin{equation}  \label{trace}
\mathrm{Trace}_{\omega_o}(i\pa\opa\varphi) \gtrsim \rho^{-2}.
\end{equation}
If the original metric $\omega_o$ does not satisfy this condition, then we replace it by $\omega_o + A \sum_{j=1}^{n+k-1}\eta_j\wedge \ol \eta_j$, where $A > 0$ is a real constant and $\eta_j,\ j=1,\ldots,n+k-1$ are $(1,0)$-forms that are smooth on $U\setminus M$ and bounded on $U$: \\
Let $\zeta_1,\ldots,\zeta_{n+k}$ be an orthonormal basis of $(1,0)$-vector fields on $U$ with respect to the original metric, smooth on $U\setminus M$ and bounded on $U$, with $\zeta_{n+k}=\xi$, where $\xi$ is as in the proof of Proposition \ref{function}. Taking $\eta_1,\ldots,\eta_{n+k}$ to be the dual basis, (\ref{trace}) is satisfied if $A$ is sufficiently big.\\

 We now modify the metric in $(\mathrm{det} N_{\mathcal{F}})^\ell$ by a bounded factor $\exp{(m\rho^2)}$. This adds to the curvature a term which is $$-mi\pa\opa\rho^2= -\sum_{j=1}^k \big( 2m\rho i\pa\opa\rho_j + 2mi\pa\rho_j\wedge\opa\rho_j\big).$$
Taking $m$ sufficiently large, and shrinking $U$ if necessary, we may therefore assume that $\mathrm{Trace}_{\omega_o}(i\Theta ((\mathrm{det} N_{\mathcal{F}})^\ell)$ negative on $U$.
Thus
$$-\mathrm{Trace}_{\omega_o}\big(i\Theta((\mathrm{det}N_{\mathcal{F}})^\ell) - ik\pa\opa\varphi\big) \gtrsim \rho^{-2}$$
on $U$.
Increasing $m$ if necessary, we may even assume that the  torsion of $\omega_o$ can be absorbed by the right-hand side of the above inequality.
But then the above estimate and the Bochner-Kodaira-Nakano inequality implies that for $u\in \mathcal{D}^{0,0}(U\setminus M, (\mathrm{det} N_{\mathcal{F}})^\ell)$ we have 
$$\int_{U\setminus M} \frac{1}{\rho^2}\vert u\vert^2 e^{k\varphi}dV_{\omega_o} \lesssim 
\int_{U\setminus M}  \vert \opa u\vert^2_{\omega_o} e^{k\varphi}dV_{\omega_o}.$$

 On the one hand, we have $\omega_o\lesssim \omega_M$ by Lemma \ref{metric}.
We even have  $dV_{\omega_M} \sim \rho^{-2k} dV_{\omega_o}$ by construction of the metric $\omega_M$.
 Therefore the preceeding inequality implies that 
$$\int_{U\setminus M} \vert u\vert^2 e^{\varphi}dV_{\omega_M} \lesssim 
\int_{U\setminus M}  \vert \opa u\vert^2_{\omega_M} e^{\varphi}dV_{\omega_M}$$
for $u\in \mathcal{D}^{0,0}(U\setminus M, (\mathrm{det} N_{\mathcal{F}})^\ell)$, which extends to $u\in L^2_{0,0}(
U\setminus M, (\mathrm{det} N_{\mathcal{F}})^\ell, \omega_M, -\varphi)$ with compact support in $U$ by the completeness of $\omega_M$ on $\ol U\setminus M$. In addition, the Levi-form of $\pa U$ has at least $n\geq 2$ negative eigenvalues. Therefore we may invoke $\cite[\mathrm{Theorem}\ 7.4]{G}$ to conclude that 

$$\Vert u\Vert_{\omega_M,-\varphi} \lesssim \Vert \opa u \Vert_{\omega_M-\varphi}$$
for all $u\in L^2_{0,0}(
U\setminus M, (\mathrm{det} N_{\mathcal{F}})^\ell, \omega_M, -\varphi)$ such that \\$\opa u\in 
L^2_{0,1}(
U\setminus M, (\mathrm{det} N_{\mathcal{F}})^\ell, \omega_M, -\varphi)$. \hfill$\square$\\

\section{Projective embeddings of tubular neighborhoods}

In \cite{O6} (see also  \cite{HM} for $k=1$), Kodaira's embedding theorem was generalized to the setting of compact Levi-flat CR manifolds, and it was shown that sufficiently high tensor powers of a positive CR line bundle over a smooth, compact Levi-flat CR manifold $M$ admit enough CR sections $s_0,\ldots, s_m$, so that the CR map $\lbrack s_0 : \ldots : s_m\rbrack$ provides a CR-embedding of $M$ into $\C\PP^m$. This applies to our situation, since $\mathrm{det} N_M$ is assumed to be positive.\\

In particular, it was proved in $\cite{O6}$ that if $\ell$ is big enough, then the $\mathcal{C}^4$-smooth CR-sections of $(\mathrm{det}N_M)^\ell$ separate the points on $M$ and give local coordinates on $M$. Using Proposition \ref{extension}, the CR-sections of $(\mathrm{det}N_M)^\ell$ can be extended to holomorphic sections of $(\mathrm{det}N_{\mathcal{F}})^\ell$ over a tubular neighborhood of $M$ in $X$.\\

Arguing by continuity, it is not difficult to see that if $\ell$ is big enough, then, after possibly shrinking $X$, the holomorphic sections of $\mathrm{det}N_M$ separate points and give local coordinates on $X$.  Hence we have  a
 holomorphic embedding $\Psi: X \hookrightarrow \C\PP^m$.\\

Using Proposition \ref{function}, we may assume that $V = \Psi(X)$ 
 is pseudoconcave in the sense of \cite{An}. Therefore we may invoke
 $\cite[\mathrm{Th\acute{e}or\grave{e}me}\ 6]{An}$ to conclude that the projective closure of $V$ in $\C\PP^m$ is a projective variety of the same dimension. By denoting $\hat X$ the projective normalization of $V$ (which exists by a theorem of Zariski), we may moreover assume that $\hat X$ is normal. $\hat X$ inherits the structure of a compact K\"ahler space from the K\"ahler metric of $\C\PP^m$. This justifies the following proposition.\\

\newtheorem{kaehler}{Proposition}[section]
\begin{kaehler} \ \\
A sufficiently small tubular neighborhood of $M$ can be holomorphically embedded into a compact (normal) K\"ahler space $\hat X$ of dimension $n+k$.
\end{kaehler}

\section{Hodge symmetry on a compact K\"ahler space}

By the results of the preceeding sections, we may assume that the Levi-flat CR manifold M is embedded into a normal compact K\"ahler space $(\hat X,\hat \omega)$ of dimension $n+k$. We define $\Omega = \mathrm{Reg}\hat X = \hat X\setminus \mathrm{Sing}\hat X$. By the normality of $\hat X$, we have $\mathrm{codim\ Sing}\hat X \geq 2$. As observed e.g. in \cite{O8}, $\Omega$ admits a complete K\"ahler metric of a particular nice form.\\

\newtheorem{metric1}{Proposition}[section]
\begin{metric1}  \label{metric1}\ \\
There exists a complete K\"ahler metric $\omega_\Omega$ on $\Omega$ and an exhaustion function $\psi: \Omega\rightarrow \lbrack 0,+\infty)$ such that
\begin{enumerate}
\item $\vert \pa\psi\vert^2_{\omega_\Omega} <\frac{1}{2(n+k)}$ 
\item The eigenvalues $\lambda_1\leq \ldots \leq \lambda_{n+k}$ of $i\pa\opa\psi$ with respect to $\omega_\Omega$ satisfy
\begin{enumerate}
\item $-1 - \frac{1}{2(n+k)} < \lambda_j < -1 + \frac{1}{2(n+k)}$ outside a compact $K_\Omega$ of $\Omega$ for $1\leq j\leq \mathrm{codim Sing}\hat X$
\item $\lambda_j < \frac{1}{2(n+k)}$ for $j > \mathrm{codim Sing} \hat X$
\end{enumerate}
\end{enumerate}
\end{metric1}

{\it Proof.} The details of the construction can be found in \cite{O8}, we only sketch the idea for the convenience of the reader.\\
Let $x_o$ be a singular point of $\hat X$, and let $f_1,\ldots, f_m$ be holomorphic functions on a neighborhood $V$ of $x_o$ that generate the ideal of holomorphic functions vanishing on $\mathrm{Sing}\hat X\cap V$ in the ring of holomorphic functions on $V$. Shrinking $V$ if necessary so that $\sum_{j=1}^m\vert f_j\vert^2 < e^{-e}$ on $V$, we set
$$\psi_V= \varepsilon \log(-\log(\sum_{j=1}^m\vert f_j\vert^2))$$
and $\omega_V = A\hat\omega -\varepsilon i\pa\opa\psi_V$ for some small $\varepsilon > 0$ and some large $A > 0$. A simple computation shows that 
$$-\varepsilon i\pa\opa\psi_V= \frac{\varepsilon i\pa\opa\log(\sum_{j=1}^m\vert f_j\vert^2)}{-\log(\sum_{j=1}^m\vert f_j\vert^2)} + \frac{i}{\varepsilon} \pa\psi_V\wedge\opa\psi_V \geq \frac{i}{\varepsilon}\pa\psi_V\wedge\opa\psi_V.$$
It is therefore easy to see that $\psi_V$ and $\omega_V$ have the required properties on $V$.\\

A detailed computation in \cite{O8} shows that $\hat\omega- i\pa\opa\psi_{V_\alpha}$ and $\hat\omega - i\pa\opa\psi_{V_\beta}$ are quasi-isometrically equivalent. Therefore the functions $\psi_\alpha$ can be patched together by a partition of unity to define a complete K\"ahler metric on $\Omega$:\\
Let $\mathcal{V} = \lbrace V_\alpha\rbrace$ be a finite open cover of $\hat X$ by such $V$, and let $\eta_\alpha$ be a smooth partition of unity associated to $\mathcal{V}$. We set
$$\psi = \underset{\alpha}{\sum} \eta_\alpha\psi_{V_\alpha}$$
and 
$$\omega_\Omega = A\hat\omega -i\pa\opa\psi.$$
As in \cite{O8}, the function and the metric thus defined have the required properties.\hfill$\square$\\

\newtheorem{finite}[metric1]{Proposition}
\begin{finite} \label{finite}\ \\
The $L^2$-Dolbeault-cohomology groups $H^{0,1}_{L^2}(\Omega,\omega_\Omega)$ and $H^{1,0}_{L^2}(\Omega,\omega_\Omega)$
are finite dimensional.
\end{finite}

{\it Proof.} The idea is to use the twisting trick of Berndtsson and Siu, since the function $\psi$ satisfies the Donnelly-Fefferman condition. On the trivial bundle $E = \Omega\times\C$ we introduce the auxiliary metric $e^{-\psi}$. Let $\lambda_1\leq\ldots\leq\lambda_{n+k}$ be the eigenvalues of $i\pa\opa\psi$ with respect to $\omega_\Omega$ as in Proposition \ref{metric1}. The curvature term $\lbrack i\pa\opa\psi,\omega_\Omega\rbrack$ acting on $(0,1)$-forms is given by $-(\lambda_2 + \ldots + \lambda_{n+k})$ (see e.g. \cite{D1}). Since $\hat X$ is normal, we have $\mathrm{codim Sing} \hat X\geq 2$. Therefore
\begin{equation} \label{curv1}
-(\lambda_2+\ldots+\lambda_{n+k}) \geq 1 -\frac{n+k-1}{2(n+k)}
\end{equation}
outside $K_\Omega$ from Proposition \ref{metric1}. Hence for every $u\in \mathcal{C}_c^{0,1}(\Omega\setminus K_\Omega)$ we have the estimate

\begin{equation}  \label{estspace1}
\big( 1-\frac{n+k-1}{2(n+k)} \big)\int_\Omega \vert u\vert^2_{\omega_\Omega} e^{-\psi} dV_{\omega_\Omega} \leq \int_\Omega ( \vert \opa u\vert^2_{\omega_\Omega} + \vert \opa^\ast_\psi u\vert^2_{\omega_\Omega}) e^{-\psi} dV_{\omega_\Omega}
\end{equation}

Now we substitute $v= u e^{-\psi/2}$. It is not difficult to see that
$$\vert \opa u\vert^2_{\omega_\Omega} e^{-\psi} \leq 2 \vert\opa v\vert^2_{\omega_\Omega} + \frac{1}{2} \vert \opa\psi\vert^2_{\omega_\Omega} \vert v\vert^2_{\omega_\Omega} \leq 2 \vert\opa v\vert^2_{\omega_\Omega} + \frac{1}{4(n+k)} \vert v\vert^2_{\omega_\Omega}.$$

Since $\opa^\ast_\psi = e^\psi \opa^\ast e^{-\psi}$, we likewise get
$$\vert\opa^\ast_\psi u\vert^2_{\omega_\Omega} e^{-\psi} \leq 2\vert\opa^\ast v\vert^2_{\omega_\Omega} + \frac{1}{2} \vert \pa\psi\vert^2_{\omega_\Omega} \vert v\vert^2_{\omega_\Omega} \leq 2\vert\opa^\ast v\vert^2_{\omega_\Omega} + \frac{1}{4(n+k)}\vert v\vert^2_{\omega_\Omega}.$$

Together with (\ref{estspace1}), these two inequalities imply
\begin{equation}  \label{estspace2}
\Vert v\Vert^2_{\omega_\Omega} \leq 4 ( \Vert\opa v\Vert^2_{\omega_\Omega} + \Vert\opa^\ast v\Vert^2_{\omega_\Omega})
\end{equation}
for all $v\in\mathcal{C}^{0,1}_c(\Omega\setminus K_\Omega)$. It is now standard to conclude that for any compact $K$ containing $K_\Omega$ in its interior, there exists a constant $C_K$ such that
$$\Vert v\Vert^2_{\omega_\Omega} \leq C_K ( \Vert \opa v\Vert^2_{\omega_\Omega} + \Vert\opa^\ast v\Vert^2_{\omega_\Omega} + \int_K \vert v\vert^2_{\omega_\Omega} dV_{\omega_\Omega})$$
for all $v\in\mathcal{C}^{0,1}_c(\Omega)$, hence for all $v\in L^2_{0,1}(\Omega,\omega_\Omega)$ by completeness of the metric $\omega_\Omega$. But the above estimate implies that $\mathcal{H}^{0,1}(\Omega,\omega_\Omega)$ is finite dimensional and isomorphic to the $L^2$-Dolbeault cohomology group $H^{0,1}_{L^2}(\Omega,\omega_\Omega)$.\\

The proof for $H^{1,0}_{L^2}(\Omega,\omega_\Omega)$ is exactly the same, since the curvature terms in degree $(0,1)$ and $(1,0)$ coincide. \hfill $\square$\\

\newtheorem{symmetry}[metric1]{Corollary}
\begin{symmetry} \label{symmetry}\ \\
We have $H^{0,1}_{L^2}(\Omega,\omega_\Omega) \simeq \overline{ H^{1,0}_{L^2}(\Omega,\omega_\Omega)}$.
\end{symmetry}

{\it Proof.} Since the metric $\omega_\Omega$ is K\"ahler, we have $\mathcal{H}^{0,1}(\Omega,\omega_\Omega) = \overline{\mathcal{H}^{1,0}(\Omega,\omega_\Omega)}$. Since $\omega_\Omega$ is complete and the $L^2$-Dolbeault cohomology groups in degrees $(0,1)$ and $(1,0)$ are finite dimensional by Proposition \ref{finite}, they are canonically isomorphic to the spaces of harmonic forms, just like in the case of compact K\"ahler manifolds. For more details, we refer the reader to \cite{D2}.\hfill $\square$\\

\section{Proof of Theorem \ref{main}}

We need a very last extension result for $\opa$-closed forms on $U$ to $\opa$-closed forms on $\Omega$ to conclude with the proof of our main result.\\

For some sufficiently small but fixed $\varepsilon > 0$ we set
$U_\varepsilon = \lbrace z\in U\mid \rho(z) < \varepsilon\rbrace$ and $V = \Omega\setminus\ol U_\varepsilon$. We may assume that $\pa V$ is smooth and that its Levi-form has $n$ positive and $k-1$ negative eigenvalues. We also set $\tilde\varphi = \tilde\max(-\log v, c\varepsilon)$, where $\tilde\max$ is a regularized maximum function and $c > 0$ is chosen such that $\tilde\varphi = \varphi$ near $\pa V$.\\

We will 
now construct a  {\it K\"ahler} metric $\omega_V$ on a neighborhood of $\ol V$ in $\Omega$, complete on $\Omega\cap\ol V$,  such that $\mu_1 + \ldots + \mu_k$ is semipositive on $V$ and $\mu_1 + \ldots + \mu_k> 1$ near $\pa V$, where $\mu_1\leq \ldots\leq \mu_{n+k}$ are the eigenvalues of $i\pa\opa \tilde\varphi$ with respect to $\omega_V$ (near $\pa V$ means that there exists some $\kappa > \varepsilon$ such that property under consideration holds on $V\setminus U_\kappa$).\\

To this extent, we define $\tilde \rho = \rho -\varepsilon$, which is a 
defining function for $\pa V$. Now we set
$$\omega_V = \omega_\Omega + \frac{1}{2\tau} i\pa\opa \tilde\rho^2$$
for some sufficiently small $\tau > 0$ ($\tilde\rho$ will be slightly modified later on). Since $\omega_V = \omega_\Omega + \frac{1}{\tau} \tilde\rho i\pa\opa\tilde\rho + \frac{1}{\tau} i\pa\tilde\rho\wedge\opa\tilde\rho$, it is clear that the condition $\mu_1´+ \ldots + \mu_k > 1$ can be satisfied on $\pa V$ for $\tau$ sufficiently small. It is also clear that $\omega_V$ is positive definite and therefore defines a K\"ahler metric on $\lbrace \tilde\rho^2 < \tau^4\rbrace$ if $\tau $ is sufficiently small.
Therefore we replace $\tilde \rho$ by $\tilde \min (\tilde\rho, a\tau^2)$, where $\tilde\min$ is a regularized minimum function. If $a>0$ is conveniently chosen, it then follows from the well known properties of regularized maximum functions, that with this new function $\tilde\rho$, the form $\omega_V$ as defined above is indeed positive definite. \\

$\omega_V$ is a K\"ahler metric, because so is $\omega_\Omega$. The assertion about the completeness follows from the completeness of $\omega_\Omega$. The assertion about$\mu_1 +\ldots +\mu_k$ is also fulfilled if the constant $c>0$ in the definition of $\tilde\varphi$ is conveniently chosen.\\

\newtheorem{harmonic}{Lemma}[section]
\begin{harmonic} \label{harmonic}\ \\
We have $\mathcal{H}_c^{0,2}(V, \omega_V, -\tilde\varphi) = \lbrace 0\rbrace \quad \big(=  L^2_{0,2}(V,\omega_V,-\tilde\varphi) \cap \mathrm{Ker}\opa_c\cap\mathrm{Ker}\vartheta\big)$.
\end{harmonic}

{\it Proof.} Taking into account that $\omega_V$ is K\"ahler, complete near $\mathrm{Sing}\hat X$ and the Levi-form of $\pa V$ has at least $n\geq 2$ positive eigenvalues, we may conclude from the Bochner-Kodaira-Nakano integral formula with boundary term (see $\cite[\mathrm{Theorem}\ 7.2]{G}$) that
\begin{equation}  \label{abc}
 \ll (\mu_1 + \ldots + \mu_{n+k-2}) u,u\gg_{\omega_V,\tilde\varphi} \leq \Vert \opa u\Vert^2_{\omega_V,\tilde\varphi} + \Vert \opa^\ast u\Vert^2_{\omega_V,\tilde\varphi} 
\end{equation}
for all $u\in L^2_{n+k,n+k-2}(V,\omega_V,\tilde\varphi) \cap\mathrm{Dom}\opa\cap\mathrm{Dom}\opa^\ast$.
Note that, since $n\geq 2$, we have $\mu_1+\ldots +\mu_{n+k-2}\geq \mu_1 +\ldots +\mu_k$. But then every harmonic form $u\in\mathcal{H}^{n+k,n+k-2}(V,\omega_V,\tilde\varphi)$ vanishes near $\pa V$. From Aronszajn's uniqueness theorem we then obtain $\mathcal{H}^{n+k,n+k-2}(V,\omega_V,\tilde\varphi)= \lbrace 0\rbrace$.
But then $\mathcal{H}_c^{0,2}(V, \omega_V, -\tilde\varphi) = \lbrace 0\rbrace$ by Serre duality (see \cite{Ch-S}).\hfill $\square$\\

\newtheorem{dolbeault}[harmonic]{Proposition}
\begin{dolbeault}  \label{dolbeault} \ \\
For every $v\in L^2_{0,2}(V,\omega_V,-\tilde\varphi)\cap\mathrm{Ker}\opa_c$ there exists
$v\in L^2_{0,1}(V,\omega_V,-\tilde\varphi)\cap\mathrm{Dom}\opa_c$ satisfying $\opa_c u = f$ on $V$.
\end{dolbeault}

{\it Proof.} If we combine the estimate (\ref{abc}) for $(n+k,n+k-1)$-forms instead of $(n+k,n+k-2)$-forms together with the twisting trick of Berndtsson and Siu as in the proof of Proposition \ref{finite}, introducing the auxiliary weight function $\psi$ near $\mathrm{Sing}\hat X$, then we obtain (as in the proof of Proposition \ref{finite}) that $H^{n+k,n+k-1}_{L^2}(V,\omega_V,\tilde\varphi)$ is finite dimensional, hence Hausdorff. It then follows by Serre duality (see \cite{Ch-S}) that $H^{0,2}_{c,L^2}(V,\omega_V,-\tilde\varphi)$ is Hausdorff. Hence
$$H^{0,2}_{c,L^2}(V,\omega_V,-\tilde\varphi) \simeq \mathcal{H}_c^{0,2}(V, \omega_V, -\tilde\varphi) = \lbrace 0\rbrace$$
by Lemma \ref{harmonic}. \hfill$\square$\\

We are now ready to prove our main result.\\

{\it Proof of Theorem \ref{main}.} 
Let $\theta = i\Theta(\mathrm{det} N_{\mathcal{F}})$ be the curvature $(1,1)$-form of the holomorphic line bundle $\mathrm{det} N_{\mathcal{F}}$ on $U$. The line bundle $\mathrm{det} N_{\mathcal{F}}$ is topologically trivial over $U$ (after possibly shrinking $U$): We cover $M$ by a finite number of balls $B_i$ such that $\mathcal{F}$ is defined by $k$ holomorphic $1$-forms with wedge product $\omega_i$ on $B_i$. Then $\mathrm{det}N_{\mathcal{F}}$ is defined by $(f_{ij})$ with $f_{ij} = \omega_i/\omega_j$. On the other hand, since $M$ is a complete intersection of $k$ real hypersurfaces $\Sigma_j$, there exist $k$ non-vanishing smooth vector fields $\xi_j$ on an open neighborhood of $M$, each pointing in the normal direction to $M$, such that they span a $k$-dimensional complex vector space. But this means that $\omega_j(\xi_1,\ldots,\xi_k)$ defines a smooth splitting of $(f_{ij})$ over an open neighborhood of $M$. Hence $\mathrm{det} N_{\mathcal{F}}$ is topologically trivial over an open neighborhood of $M$.\\

But then there exists a smooth real-valued 1-form $\beta$ on $U$ such that
$$d\beta = \theta = \pa\beta^{0,1} + \opa\beta^{1,0},$$
where $\beta = \beta^{1,0} + \beta^{0,1}$ can be chosen to be real ($\beta^{0,1} = \overline{\beta^{1,0}}$). From standard type considerations we get $\opa\beta^{0,1} = 0$ on $U$. We choose a smooth cut-off function $\chi$ which is equal to one on $U\setminus V$ and compactly supported in $U$. Then $f = \opa(\chi\beta^{0,1})$ is $\opa_c$-closed and belongs to $L^2_{0,2}(V, \omega_V, -\tilde\varphi)$. From Proposition \ref{dolbeault}
we obtain $u\in L^2_{0,1}(V, \omega_V, -\tilde\varphi)$ satisfying $\opa_c u = f$ on $V$. But then the extension of $u$ to $\Omega$, still denoted by $u$, obtained by setting it zero outside $V$ satisfies $\opa u = f$ on $\Omega$ in the distribution sense. Hence $\tilde\beta^{0,1} = \chi\beta^{0,1} - u \in L^2_{0,1}(\Omega,\omega_\Omega)$ satisfies $\opa\tilde\beta^{0,1}=0$ in $\Omega$
 But then Corollary \ref{symmetry} implies that there exist forms $\alpha\in L^2_{0,0}(\Omega,\omega_\Omega), \eta\in L^2_{0,1}(\Omega,\omega_\Omega)$, $\pa\eta = 0$, such that $\tilde\beta^{0,1} = \eta + \opa\alpha$ on $\Omega$. Moreover, we can without loss of generality assume that $\alpha$ is smooth on $U$ (since $\beta^{0,1}$ is smooth there)
Therefore, setting $\phi = i(\overline{\alpha} -\alpha)$, one obtains $\theta = i\pa\opa\phi$ on every leaf of the Levi-foliation of $M$. The existence of a potential for the positive curvature is, however, a contradiction to the maximum principle on the leaves of the foliation. \hfill $\square$\\

{\bf Acknowledgements.} 
The research on this project was supported by Deutsche Forschungsgemeinschaft (DFG, German Research Foundation, grant BR 3363/2-2).
I would also like to thank the anonymous referee for the careful reading and for suggesting several improvements.\\

\end{document}